\documentclass[conference]{IEEEtran}
\IEEEoverridecommandlockouts
\usepackage{cite}
\usepackage{amsmath,amssymb,amsfonts, bm}
\usepackage{algorithmic}
\usepackage{graphicx}
\usepackage{textcomp}
\usepackage{xcolor}
\usepackage{stfloats}
\usepackage{algorithm}
\usepackage{graphicx}
\usepackage{subfigure}

\def\BibTeX{{\rm B\kern-.05em{\sc i\kern-.025em b}\kern-.08em
    T\kern-.1667em\lower.7ex\hbox{E}\kern-.125emX}}
\begin{document}

\title{A Collaborative Jamming Algorithm Based on Multi-UAV Scheduling
}

\author{
	\IEEEauthorblockN{Yixin Jiang\IEEEauthorrefmark{1}, Lingyun Zhou\IEEEauthorrefmark{1},  Yijia Tang\IEEEauthorrefmark{1}, Ya Tu\IEEEauthorrefmark{1}, Chunhong Liu\IEEEauthorrefmark{2}, Qingjiang Shi\IEEEauthorrefmark{1}\IEEEauthorrefmark{3}} 
	\IEEEauthorblockA{\IEEEauthorrefmark{1}School of Software Engineering, Tongji University, Shanghai 201804, China} 
		\IEEEauthorblockA{\IEEEauthorrefmark{2}Science and Technology on Communication Information Security Control Laboratory, Jiaxing 314033, China}
	\IEEEauthorblockA{\IEEEauthorrefmark{3}Shenzhen Research Institute of Big Data, Shenzhen 518172, China}  
}

\maketitle

\begin{abstract}
In this paper, we consider the problem of multi-unmanned aerial vehicles' (UAVs)' scheduling for cooperative jamming, where UAVs equipped with directional antennas perform collaborative jamming tasks against several targets of interest. To ensure effective jamming towards the targets, we formulate it as an non-convex optimization problem, aiming to minimize the communication performance of the targets by jointly optimizing UAVs’ deployment and directional antenna orientations. 
Due to the unique structure of the problem, we derive an equivalent transformation by introducing a set of auxiliary matrices. Subsequently, we propose an efficient iterative algorithm based on the alternating direction method of multipliers (ADMM), which decomposes the problem into multiple tractable subproblems solved in closed-form or by gradient projection method.
Extensive simulations validate the efficacy of the proposed algorithm.
\end{abstract}


\section{INTRODUCTION}

Unmanned aerial vehicles (UAVs), being large-scale collaborative units, have many  civilian and military applications due to advantages such as flexible deployment, adaptable coverage, and line-of-sight (LoS) connectivity \cite{UAV_application_v1,UAV_application_v2}. 
The effective collaboration among multiple UAVs not only compensates for the limited capability of a single unit, but also enhances the fault tolerance of the entire system, thereby ensuring effectiveness and feasibility in completing various tasks.
Due to these satisfactory characteristics, collaborative operations based on UAV swarm have recently attracted great interest in a wide range of fields \cite{UAVs_colla_v1,UAVs_colla_v2}.

The future battlefield environment, characterized by high
mobility, intense adversarial activities, and increased informatization, has put forward various collaborative applications of UAVs, in which jamming hostile targets is of paramount importance \cite{UAV_Intere}. 
However, due to the adversarial environment and responsiveness requirement, the design of an effective target interference strategy is a formidable challenge. 
Thus, it is imperative to devise efficient resource allocation strategy, thereby enabling more precise and powerful interference. 

Recently, numerous studies have focused on jamming strategies applied by multi-UAV systems \cite{UAV_IS_v1, UAV_IS_v2, UAV_IS_v3}. 
Specifically, to achieve efficient interference with communication signals of the targets, an energy-aware tracking and jamming framework was proposed to jointly optimize the mobility and jamming power of the multiple UAVs \cite{UAV_IS_v1}.
In \cite{UAV_IS_v2}, a UAV detector and jammer system was explored, in which a directional jammer transmits a powerful signal to disrupt the communication between a malicious UAV and its controller upon identifying the UAV.
Attempts to investigate an innovative and cost-effective autonomous anti-UAV system designed to detect, track, and jam targets was made in \cite{UAV_IS_v3}.

Although the previous studies have made valuable contributions to interference methods with multi-UAV scheduling, the design of an efficient, reliable and powerful jamming algorithm remains relatively unexplored in the existing literature.
Motivated by this, our paper delves into the coordinated resource allocation for precise and effective interference in multi-UAV system. 
Specifically, we consider a multi-UAV interference scenario, wherein each UAV is equipped with a directional antenna for collaborative disruption and interference against multiple targets.
The main challenge of this problem is how to coordinate multi-UAVs' spatial resources to simultaneously jamming multiple targets of interest.
We formulate a non-convex optimization problem, aiming to minimize the average signal-to-interference-plus-noise ratio (SINR) of the targets by jointly optimizing UAVs' deployment and directional antenna orientations under stringent constraints on resource availability. 
To effectively tackle the above problem, we first introduce variables to decouple the non-convex constraints and reformulate it. Then, by virtue of the particular problem structure, we adopt an alternating direction multipliers method (ADMM)-based method for the problem. The key component of the algorithm is to update UAVs' deployment, in which we propose and adopt an improved gradient projection method.
The simulations clearly show the enhanced interference effectiveness of the proposed scheme and deliver useful insights for practical system design.

The remaining sections of this paper are structured as follows. Section II outlines the interference scenario based on UAV swarms and presents the formulated optimization problem. Section III introduces an iterative algorithm with low complexity for the identified problem. The results of our proposed algorithm analyzed through simulations are presented in Section IV. Finally, section V concludes the paper.

\begin{figure}[!t]
	\centering
	\includegraphics[width=8.6cm,height=2.7in]{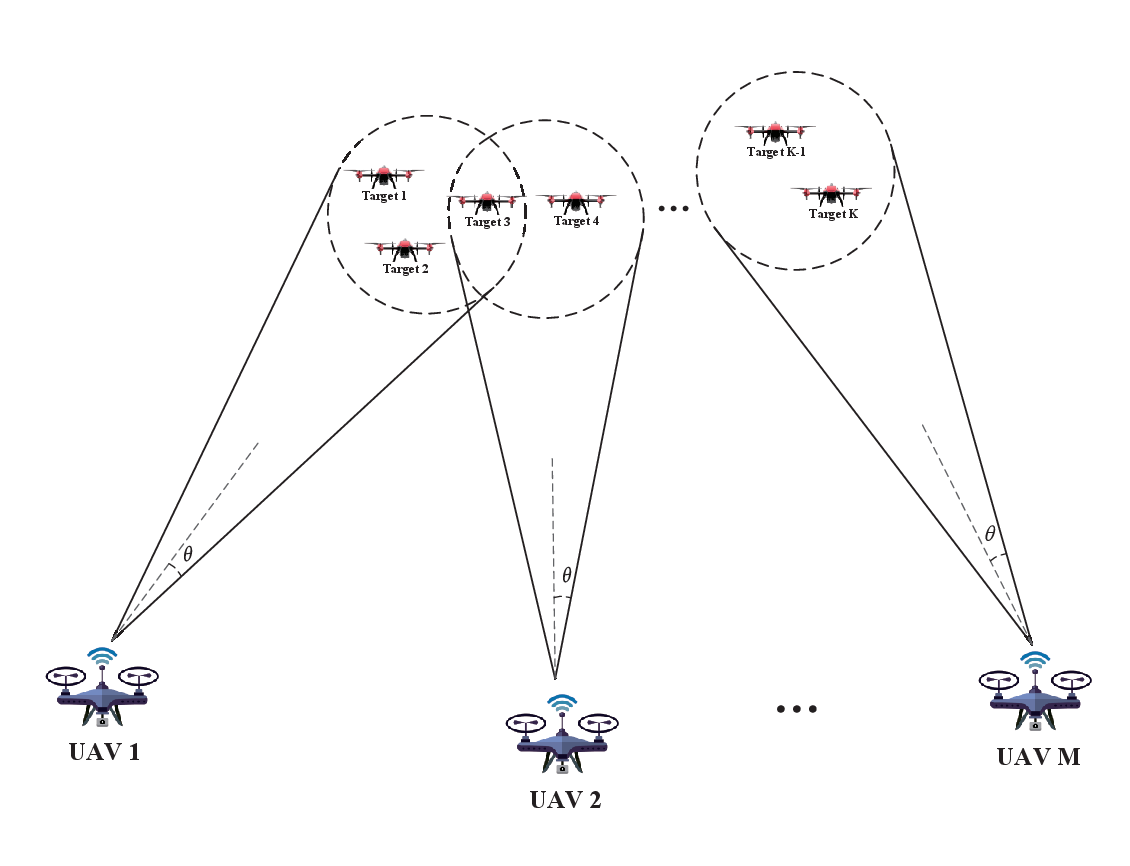}
	\caption{Illustration of a UAV swarm-based interference network scenario, where $M$ UAV cooperatively perform interference to $K$ targets.} 
	\label{fig:UAV_Scenario}
\end{figure}

\section{SYSTEM MODEL AND PROBLEM FORMULATION}
As shown in Fig. \ref{fig:UAV_Scenario}, we consider a multi-UAV jamming network, where each UAV is equipped with a directional antenna to perform collaborative interference against opposing targets. Specifically, we adopt a 3-dimensional (3D) Cartesian coordinate framework encompassing $M$ UAVs and $K$ targets, denoted by $i \in \mathcal{M} \triangleq\{1, \cdots, M\}$ and $k \in \mathcal{K} \triangleq\{1, \cdots, K\}$ respectively.
Each UAV maintains a consistent altitude $H$ to avoid frequent ascents and descents and minimize energy consumption.
As such, let $\boldsymbol{Q}\triangleq\left[\boldsymbol{q}_{1}, \boldsymbol{q}_{2}, \ldots, \boldsymbol{q}_{M}\right]^T\in \mathbb{R}^{M \times 3}$ denote the deployment matrix  with
$\boldsymbol{q}_{i} = \left[x_{i}, y_{i}, H\right]$ being the coordinates of the $i$-th UAV.
Similarly, we define $\boldsymbol{Q}_t\triangleq\left[\boldsymbol{q}_{t,1}, \boldsymbol{q}_{t,2}, \ldots, \boldsymbol{q}_{t,K}\right]^T\in \mathbb{R}^{M \times 3}$ to represent the location matrix of targets with $\boldsymbol{q}_{t, k} = \left[x_{t, k}, y_{t, k}, z_{t, k}\right]$ as the location of target $k$. 


In the examined multi-UAV cooperative interference scenario, the directional antennas mounted on the UAVs play a pivotal role in bolstering the channel gain in specific directions. Typically, the directional antenna gain within the main lobe range $(-\theta, \theta)$ is estimated to hold a particular value, while beyond this range, it approaches zero or a negligible quantity. This characteristic gives rise to an approximate directional radiation gain as:
\begin{equation}
	\label{radiation_gain}
	G_{i,k} =  \begin{cases}10 \log _{10}\left(\exp \left(-\frac{\alpha_{i,k}^2}{2 \theta^2}\right)\right) &\alpha_{i,k}\in(-\theta, \theta)  \\ 0, &otherwise \end{cases}
\end{equation}
\noindent where ${\alpha}_{i,k} =  arccos(cos(\phi_{i, k}-\psi_i)cos(\varphi_{i, k}))$, with $\phi_{i,k}$ and $\varphi_{i,k}$ denoting the horizontal angle and pitch angle from UAV $i$ to target $k$, respectively. 

Specifically, $\phi_{i,k}=\arctan \left(x_{t, k}-x_{i}, y_{t, k}-y_{i}\right)$, and $\varphi_{i,k}=\arctan \big(z_{t, k}-H, \sqrt{\left(y_{t, k}-y_{i}\right)^2 +\left(x_{t, k}-x_{i}\right)^2}\big)$. Note that in our analyzed scenario, the directional antenna of each UAV is fixed in the vertical direction. Therefore, the  angle-dependent radiation gain remains unaffected by the antenna's elevation angle and solely relies on the azimuth angle. 

Assuming that all targets receive signals from the control center located at $\boldsymbol{q}_s \triangleq\left[x_s, y_s, z_s\right]$, and the transmission power is $P_s$.
Here, the channel gain of the LoS links between the target and the control center can be characterized by the free-space channel model, as expressed below:
\begin{equation}
h_{s, k}=\beta_0 d_{s, k}^{-2}=\frac{\beta_0}{\left\|\boldsymbol{q}_{t, k}-\boldsymbol{q}_s\right\|_2^2}.
\end{equation}

Hence, representing the interference power transmitted by UAV $i$ as $P_i$, the received SINR for target $k$ can be distinctly characterized as:
\begin{equation}
	\gamma_k(\boldsymbol{Q},\boldsymbol{\psi})=\frac{P_s h_{s, k}}{\sum_{i=1}^M \beta_0 P_i N_i G_{i, k}\left(\boldsymbol{q}_i, \psi_i\right) d_{i,k}^{-2}(\boldsymbol{q}_i)+\sigma_k^2},
\end{equation}
where $\sigma_k^2$ denotes the power of AWGN at the receiver side, $\beta_{0}$ represents the channel power gain at a specific unit distance, $N_i$ is the count of isotropic antenna elements, and $d_{i,k}\left(\boldsymbol{q}_{i}\right)=\left\|\boldsymbol{q}_{t,k}-\boldsymbol{q}_{i}\right\|_{2}$ indicates the distance from target $k$ to UAV $i$. 

To evaluate the jamming performance of the entire UAV swarm-based network, introduce the average SINR function of the targets as:
\begin{equation}
	\gamma(\boldsymbol{Q}, \boldsymbol{\psi})=\frac{1}{K}\sum_{k=1}^K \gamma_k(\boldsymbol{Q}, \boldsymbol{\psi}).
\end{equation}

It is evident that the lower $\gamma(\boldsymbol{Q}, \boldsymbol{\psi})$ is, the better the jamming performance of the UAV swarm-based network achieves. 
Motivated by this, the goal of this paper is to minimize the average SINR function of the targets by jointly optimizing the placement locations (i.e., $\boldsymbol{Q}$), along with the directional antenna orientations (i.e., $\boldsymbol{\psi}$). Thus, the optimization problem can be formulated as:
\begin{subequations}  \label{eq:OriFor}
	\begin{align}
		\underset{\boldsymbol{Q}, \boldsymbol{\psi}}{\operatorname{min}} \quad & 
		\gamma \left(\boldsymbol{Q}, \boldsymbol{\psi} \right) \\
		\text { s.t. } \quad 
		& \boldsymbol{q}_{i} \in \mathcal{D}_{q}, \quad \forall i, \label{subeq:OriFor_b} \\
		& \left\|\boldsymbol{q}_{i}-{\boldsymbol{q}}_{t,k}\right\|_{2} \geqslant S_{l}, \quad \forall i, k, \label{subeq:OriFor_c}\\
		& \left\|\boldsymbol{q}_{i}-\boldsymbol{q}_{j}\right\|_{2} \geqslant R_{l}, \quad \forall i \neq j, \label{subeq:OriFor_d} \\
		& -\pi < \psi_{i} \leq \pi, \quad \forall i, \label{subeq:OriFor_f}
	\end{align}
\end{subequations}
where \eqref{subeq:OriFor_b} represents the deployable area $ \mathcal{D}_{q}$.
For the sake of effectively interfering with targets while mitigating vulnerabilities, \eqref{subeq:OriFor_c} defines the minimum separation distance (i.e., $S_{l}$) between targets and UAVs. 
Furthermore, \eqref{subeq:OriFor_d} establishes the anti-collision boundary with $R_{l}$, the minimum distance between any two UAVs.
\eqref{subeq:OriFor_f} denotes directional antenna orientation constraint of each UAV.

\section{PROPOSED ALGORITHM}
Owing to the numerous mutually coupled non-convex constraints in the problem, conventional optimization algorithms struggle to seek straightforward solutions.
ADMM offers advantages such as rapid convergence and commendable convergence performance \cite{admm_convergence}, and excels in tackling such non-convex problems of this nature. The formulation of the problem and the solving process is elaborated below.

Concerning the separation distance constraint between the interference UAVs and the targets in the model, we articulate it as follows:
\begin{equation}
	\boldsymbol{b}_{i, k}=\boldsymbol{q}_i-\boldsymbol{q}_{t, k}, \forall i, k.
\end{equation}

Define $\boldsymbol{B}  \in \mathbb{R}^{M*K \times 3}$ as the matrix composed by all $\boldsymbol{b}_{i,k}$ elements, its feasible domain is $\mathcal{D}_{b}=\left\{\boldsymbol{b}_{i, k} \mid\left\|\boldsymbol{b}_{i, k}\right\|_2\geq S_l\right\}$. The constraint can be equally written in the following manner:
\begin{equation}
	 \boldsymbol{A}_{1}\boldsymbol{Q}-\boldsymbol{A}_{2}\boldsymbol{Q}_{t}=\boldsymbol{B},
\end{equation}
where $\boldsymbol{A}_{1}\in \mathbb{R}^{M*K \times M}$,  $\boldsymbol{A}_{2}\in \mathbb{R}^{M*K \times K}$ denotes the coefficient matrix in the reformulation process.

Next, we reformulate the anti-collision constraint between UAVs. Define variables:
\begin{equation}
		\boldsymbol{c}_{i, j}=\boldsymbol{q}_i-\boldsymbol{q}_j, \forall i \neq j.
\end{equation}

Define $\boldsymbol{C} \in \mathbb{R}^{\frac{M(M-1)}{2} \times 3}$ as the matrix composed by all $\boldsymbol{c}_{i,j}$ elements, its feasible domain is $\mathcal{D}_{c}=\left\{\boldsymbol{c}_{i, j} \mid\left\|\boldsymbol{c}_{i, j}\right\|_2\geq R_l\right\}$. The constraint is expressed through the subsequent equation:
\begin{equation}
	\boldsymbol{A}_{3}\boldsymbol{Q}=\boldsymbol{C},
\end{equation}
where $\boldsymbol{A}_{3}\in \mathbb{R}^{\frac{M(M-1)}{2} \times M}$ denotes the coefficient matrix in the reformulation process.

Consequently, the original optimization problem \eqref{eq:OriFor} can be equivalently transformed as follows:
\begin{subequations}  \label{eq:ProFor}
	\begin{align}
		\underset{\boldsymbol{Q}, \boldsymbol{\psi}}{\operatorname{min}} \quad & 
		\gamma \left(\boldsymbol{Q}, \boldsymbol{\psi} \right) \\
		\text { s.t. } \quad 
		& \boldsymbol{A}_{1}\boldsymbol{Q}-\boldsymbol{A}_{2}\boldsymbol{Q}_{t}=\boldsymbol{B}, \label{subeq:ProFor_b} \\
		& \boldsymbol{A}_{3}\boldsymbol{Q}=\boldsymbol{C}, \label{subeq:ProFor_c}\\
		&\boldsymbol{q}_{i} \in \mathcal{D}_{q}, \quad \forall i, \label{subeq:ProFor_d} \\
		&\boldsymbol{B} \in \mathcal{D}_{b}, \label{subeq:ProFor_e} \\
		&\boldsymbol{C} \in \mathcal{D}_{c}, \label{subeq:ProFor_f} \\
		& -\pi \leq \psi_{i} \leq \pi, \quad \forall i. \label{subeq:ProFor_g}
	\end{align}
\end{subequations}

For the above problem, the scaled form of augmented Lagrangian function can be written as:
\begin{equation}
	\begin{aligned}
	 &\quad L_\rho\left(\boldsymbol{Q}, \boldsymbol{\psi}, \boldsymbol{B}, \boldsymbol{C} ; \boldsymbol{\chi}, \boldsymbol{\mu}, \rho_1, \rho_2\right)\\
	 	&\quad =\gamma(\boldsymbol{Q}, \boldsymbol{\psi})+\frac{\rho_1}{2}\left\|\mathbf{A}_1 \boldsymbol{Q}-\mathbf{A}_2\boldsymbol{Q}_t-\boldsymbol{B}+\boldsymbol{\chi}\right\|_2^2- 
		 \frac{\rho_1}{2}\|\boldsymbol{\chi}\|_2^2\\
		 &\quad +\frac{\rho_2}{2}\left\|\mathbf{A}_3 \boldsymbol{Q}-\boldsymbol{C}+\boldsymbol{\mu}\right\|_2^2-\frac{\rho_2}{2}\|\boldsymbol{\mu}\|_2^2,
	\end{aligned}
\end{equation}
wherein $\boldsymbol{\chi}$, $\boldsymbol{\mu}$ are Lagrangian multiplier matrices, $\rho_1$, $\rho_2$ are penalty factors. Define $\boldsymbol{Q}_l$, $\boldsymbol{\psi}_l$, $\boldsymbol{B}_l$, $\boldsymbol{C}_l$, $\boldsymbol{\chi}_l$, $\boldsymbol{\mu}_l$ to signify values of parameters in the $l$-th iteration respectively. 
Drawing from the previously outlined definitions, the subsequent section elaborates on the ADMM-based variable update procedure during the $l$-th iteration.

Step $1$: Update variable $\boldsymbol{B}$, $\boldsymbol{C}$, that is, solve the following problems simultaneously in parallel:
\begin{equation}\label{resolve_b}
  \begin{aligned}
	 & \min _{\boldsymbol{B}}\left\|\mathbf{A}_1 \boldsymbol{Q}^l-\mathbf{A}_2\boldsymbol{Q}_t-\boldsymbol{B}+\boldsymbol{\chi}^l\right\|_2^2 \\ & \text { s.t. } \boldsymbol{B} \in \mathcal{D}_{b} . 	 
	 \end{aligned}
 \end{equation}
 \begin{equation}
 	\begin{aligned}\label{resolve_c}
 	  & \min _{\boldsymbol{C}}\left\|\mathbf{A}_3 \boldsymbol{Q}^l-\boldsymbol{C}+\boldsymbol{\mu}^l\right\|_2^2 \\ & \text { s.t. } \boldsymbol{C} \in \mathcal{D}_{c}.	 
 	\end{aligned}
 \end{equation}
 
	Owing to the structure of the variable $\boldsymbol{B}$, problem \eqref{resolve_b} can be broken down into $M * K$ subproblems for parallel solving. Among them, the $v$-th subproblem is as follows:
	\begin{equation}
			\begin{aligned} & \min _{\boldsymbol{b}_v}\left\|\boldsymbol{q}_{j 1}^l-\boldsymbol{q}_{t, j 2}-\boldsymbol{b}_v+\boldsymbol{\chi}_v^l\right\|_2^2 \\ & \text { s.t. }\left\|\boldsymbol{b}_v\right\|_2 \geq S_l .\end{aligned}
	\end{equation}
	
	The closed-form solution for the above problem can be expressed as:
	\begin{equation}
			\boldsymbol{b}_v^{l+1}=\left\{\begin{array}{c}\frac{\max \left(\xi_1, S_l\right)}{\xi_1}\left(\boldsymbol{q}_{j 1}^l-\boldsymbol{q}_{t, j 2}+\boldsymbol{\chi}_v^l\right) \\ \xi_1=\left\|\boldsymbol{q}_{j 1}^l-\boldsymbol{q}_{t, j 2}+\boldsymbol{\chi}_v^l\right\|_2\end{array}\right.
	\end{equation}
	
Similarly, problem \eqref{resolve_c} can be split into $\frac{M(M-1)}{2} $ subproblems, and the $v$-th is:
\begin{equation}
\begin{aligned} & \min _{\boldsymbol{c}}\left\|\boldsymbol{q}_{j 1}^l-\boldsymbol{q}_{t, j 2}-\boldsymbol{c}_v+\mu_v^l\right\|_2^2 \\ & \text { s.t. }\left\|\boldsymbol{c}_v\right\|_2 \geq R_l .\end{aligned}
\end{equation}

	The closed-form solution for the above problem can be expressed as:
	\begin{equation}
		\boldsymbol{c}_v^{l+1}=\left\{\begin{array}{c}\frac{\max \left(\xi_2, S_l\right)}{\xi_2}\left(\boldsymbol{q}_{j 1}^l-\boldsymbol{q}_{j 2}^l+\mu_v^l\right) \\ \xi_2=\left\|\boldsymbol{q}_{j 1}^l-\boldsymbol{q}_{j 2}^l+\mu_v^l\right\|_2\end{array}\right.
	\end{equation}
	
Step $2$: Update $\boldsymbol{Q}$, equivalent to solving the following problem:
\begin{equation}\label{resolve_q}
\begin{aligned} & \min _{\boldsymbol{Q}} \gamma\left(\boldsymbol{Q}, \boldsymbol{\psi}^l\right)+\frac{\rho_1}{2}\left\|\mathbf{A}_1 \boldsymbol{Q}-\mathbf{A}_2\boldsymbol{Q}_t-\boldsymbol{B}^l+\boldsymbol{\chi}^l\right\|_2^2-\frac{\rho_1}{2}\left\|\boldsymbol{\chi}^l\right\|_2^2 \\
	&+\frac{\rho_2}{2}\left\|\mathbf{A}_3 \boldsymbol{Q}-\boldsymbol{C}^l+\boldsymbol{\mu}^l\right\|_2^2  -\frac{\rho_2}{2}\left\|\boldsymbol{\mu}^l\right\|_2^2 \\
	& \text { s.t. } \boldsymbol{q}_i \in \mathcal{D}_{q}.   \end{aligned}
\end{equation}

After applying appropriate approximation transformation to $G_{i,k}$ in \eqref{radiation_gain}, the objective function becomes continuous and differentiable with respect to $\boldsymbol{Q}$.  
Denoting the objective function of the problem as $P_{\text{admm}}(\boldsymbol{Q})$, the gradient projection algorithm can be utilized for tackling the problem. Due to the non-convexity and complexity of the problem, original algorithm takes a significant amount of time to solve. This is mainly caused by inappropriate step lengths, suboptimal descent directions and oscillation. Improved gradient descend algorithms, widely applied in the field of machine learning, can tackle the aforementioned problems to a great extent\cite{Two_gradient,NAG}. Although they are originally devised for unconstrained optimization problems, in real-world applications the constraint space can often be approximated as a plane. By choosing the coordinate system appropriately, we can simplify the projection onto the constrained space to limit UAVs' x coordinates. We integrate Nesterov Acceleration Gradient(NAG) and Root Mean Square Propagation(RMSprop) into the algorithm and make such an adaptive adjustment: Due to the simplified UAV deployment constraint, it's essential to optimize the x and y coordinates separately. The procedure is concluded as follows:
First, we compute the gradient $\boldsymbol{G}$ of $P_{\text{admm}}(\boldsymbol{Q})$ with respect to $\boldsymbol{Q}$;
Second, calculate the hybrid gradient descent separately in the x direction and y direction;
 Then, calculate the projection of $\overline{\boldsymbol{Q}}$ onto the constraint space $\boldsymbol{\Omega}_{D_{\mathrm{q}}}$, denoted as $\boldsymbol{Q}_{proj}$;
 Last, choose a suitable step size $\alpha$ and search in the direction of $\boldsymbol{Q}_{proj}-\boldsymbol{Q}$ to update the variable $\boldsymbol{Q}$.
 Repeat the above steps until the algorithm converges. The detailed procedure are concluded in Algorithm 1. 
 \begin{algorithm}[!h]
 	\caption{Gradient Projection for Problem \eqref{resolve_q}}
 	\label{gradient_projection}
 	\begin{algorithmic}[1]
 		\STATE \textbf{Initialize} $\boldsymbol{Q}$ , $\boldsymbol{D}_\text{last} = 0$, $\boldsymbol{G}_\text{last} = 0$, $v_\text{x,last} = 0$,$v_\text{y,last} = 0$, $t = 1$. Define the maximum iteration number $T_{1,\max}$, and set threshold for iteration termination.
 		\REPEAT
 		\STATE Compute gradient:
 		$\boldsymbol{G}=\nabla_{\boldsymbol{Q}}\left[P_\text{admm }(\boldsymbol{Q})\right] ;$
 	\STATE Calculate hybrid descent:
 	
 	$\begin{aligned}
 			& \boldsymbol{D}=\beta_\text{NAG} \boldsymbol{D}_\text{last} +\boldsymbol{G}+\beta_\text{NAG}\left(\boldsymbol{G}-\boldsymbol{G}_ \text{last }\right)		
 				\\ &
 			v_\text{x}=\frac{\rho_\text{RMSProp} v_\text{x,last}+(1-\rho_\text{RMSProp}) \left\|\boldsymbol{G}_\text{x}\right\|_2^2}{1-\rho_\text{RMSProp}^t}  			
 		\\ &
 		v_\text{y}=\frac{\rho_\text{RMSProp} v_\text{y,last}+(1-\rho_\text{RMSProp}) \left\|\boldsymbol{G}_\text{y}\right\|_2^2}{1-\rho_\text{RMSProp}^t} \\ &
 		\overline{\boldsymbol{Q}}_\text{x}=\boldsymbol{Q}_\text{x}-\frac{\alpha_\text{NAG} \boldsymbol{D}_\text{x}}{\sqrt{v}_\text{x}+\epsilon_\text{RMSProp}}
 		 \\ &
 		\overline{\boldsymbol{Q}}_\text{y}=\boldsymbol{Q}_\text{y}-\frac{\alpha_\text{NAG} \boldsymbol{D}_\text{y}}{\sqrt{v_\text{y}}+\epsilon_\text{RMSProp}}
 	\end{aligned}$
 		\STATE Calculate projection onto constraint space:
 		
 		$\boldsymbol{Q}_\text{proj}=P_{{\Omega}_{D}}(\overline{\boldsymbol{Q}}) ;$
 		\STATE Update $\boldsymbol{Q}$ according to following strategy:
 		
 		$\boldsymbol{Q} \leftarrow-\boldsymbol{Q}+\alpha(\boldsymbol{Q}_\text{proj}-\boldsymbol{Q}) ;$
 		\STATE Update variables:
 		$\boldsymbol{D}_\text{last}\leftarrow \boldsymbol{D};$
 		
 			$\boldsymbol{G}_\text{last}\leftarrow \boldsymbol{G};$
 			
 				$v_\text{x,last}\leftarrow v_\text{x};$
 				
 					$v_\text{y,last}\leftarrow v_\text{y};$
 				
 				$t_\text{last} \leftarrow t;$
 		\UNTIL{the objective function converges, or the maximum iteration number is reached.}		
 	\end{algorithmic}
 \end{algorithm}
 
 Step 3: Update variable $\boldsymbol{\psi}$, that is, solving the following problem:
 \begin{equation}
 	 \begin{aligned} & \min _{\boldsymbol{\psi}}\gamma\left(\boldsymbol{Q}^l, \boldsymbol{\psi}\right) \\ & \text { s.t. }-\pi<\psi_i \leq \pi, \forall i,\end{aligned}
 \end{equation}
 
In addressing this problem, the periodicity of antenna angles within the objective function leads us to bypass the constraints and simply opt for a gradient descent algorithm with varied initial values. Afterward, the outcomes is converted into the constraint interval.
 
 Step 4: Update Lagrangian multipliers $\boldsymbol{\chi}, \boldsymbol{\mu}$ as follows:
 \begin{equation}
 	 \boldsymbol{\chi}^{l+1}=\left\{\begin{array}{c}\tilde{\boldsymbol{\chi}}^{l+1}, \text { if } \chi_{\max }^{l+1}<\varpi_\chi \\ \tilde{\boldsymbol{\chi}}^{l+1} / \chi_{\max }^{l+1}, \text { else }\end{array}\right.
 \end{equation}
\begin{equation}
\boldsymbol{\mu}^{l+1}=\left\{\begin{array}{c}\tilde{\boldsymbol{\mu}}^{l+1}, \text { if } \mu_{\max }^{l+1}<\varpi_\mu \\ \tilde{\boldsymbol{\mu}}^{l+1} / \mu_{\max }^{l+1}, \text { else }\end{array}\right.
\end{equation}
In the expression, $\tilde{\boldsymbol{\chi}}^{l+1}=\boldsymbol{\chi}^l+\mathbf{A}_1 \boldsymbol{Q}^{l+1}-\mathbf{A}_2\boldsymbol{Q}_t-\boldsymbol{B}^{l+1}, \tilde{\boldsymbol{\mu}}^{l+1}=\boldsymbol{\mu}^l+\mathbf{A}_3 \boldsymbol{Q}^{l+1}-\boldsymbol{C}^{l+1}$. $\varpi_\mu$ and $\varpi_\chi$ are positive numbers that are big enough.

Stop condition: Continue the aforementioned steps, until the following condition is met:
\begin{small}
	\begin{equation}
		\varepsilon^{l+1}=\left\|\boldsymbol{A}_1 \boldsymbol{Q}^{l+1}-\boldsymbol{A}_2\boldsymbol{Q}_t-\boldsymbol{B}^{l+1}\right\|_2+\left\|\boldsymbol{A}_3 \boldsymbol{Q}^{l+1}-\boldsymbol{C}^{l+1}\right\|_2 \leq \eta
	\end{equation}
\end{small}
\section{NUMERICAL RESULTS}

\begin{figure*}[!h]
	\centering
	\subfigure[$\text{ Single UAV interferes with 3 targets.}$]{
		\label{fig:subfig:a}
		\includegraphics[width=8.8cm,height=2.5in]{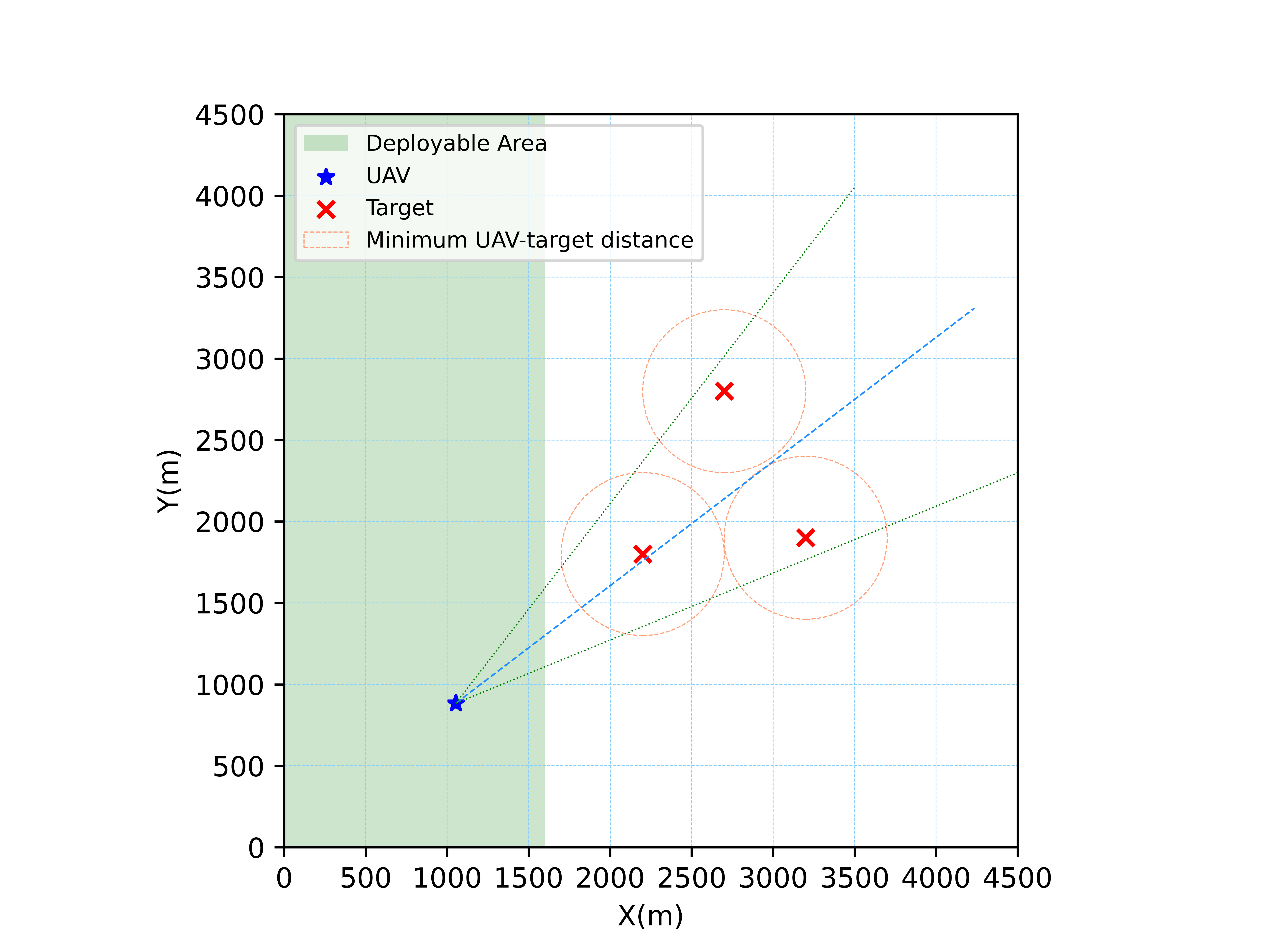}}
	\hspace{0.01in}
	\subfigure[$\text{ 2 UAVs cooperatively interfere with 3 targets.}$]{
		\label{fig:subfig:b}
		\includegraphics[width=8.8cm,height=2.5in]{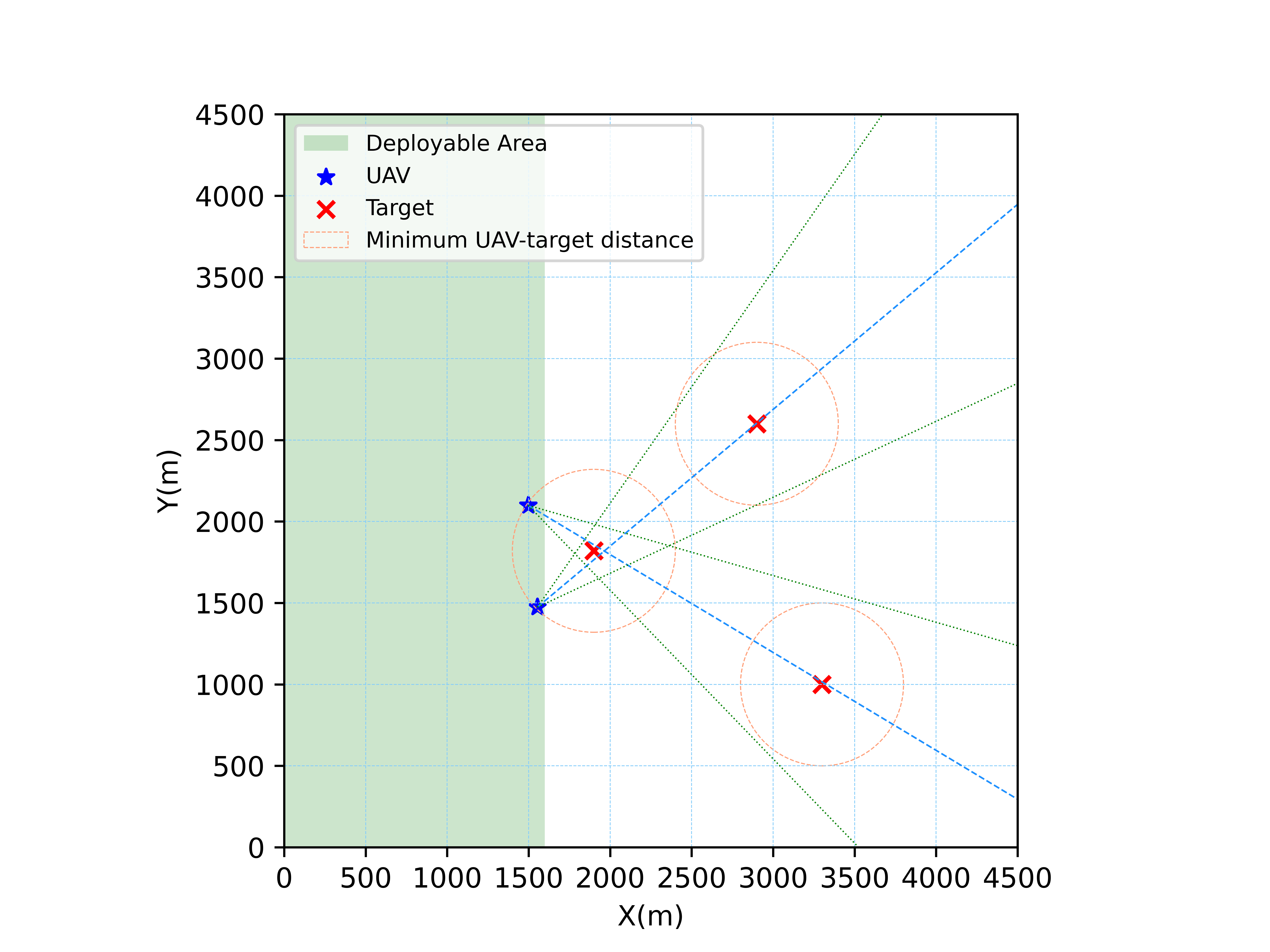}}
	\hspace{0.01in}
	\subfigure[$\text{ 3 UAVs cooperatively interfere with 3 targets.}$]{
		\label{fig:subfig:c}
		\includegraphics[width=8.8cm,height=2.5in]{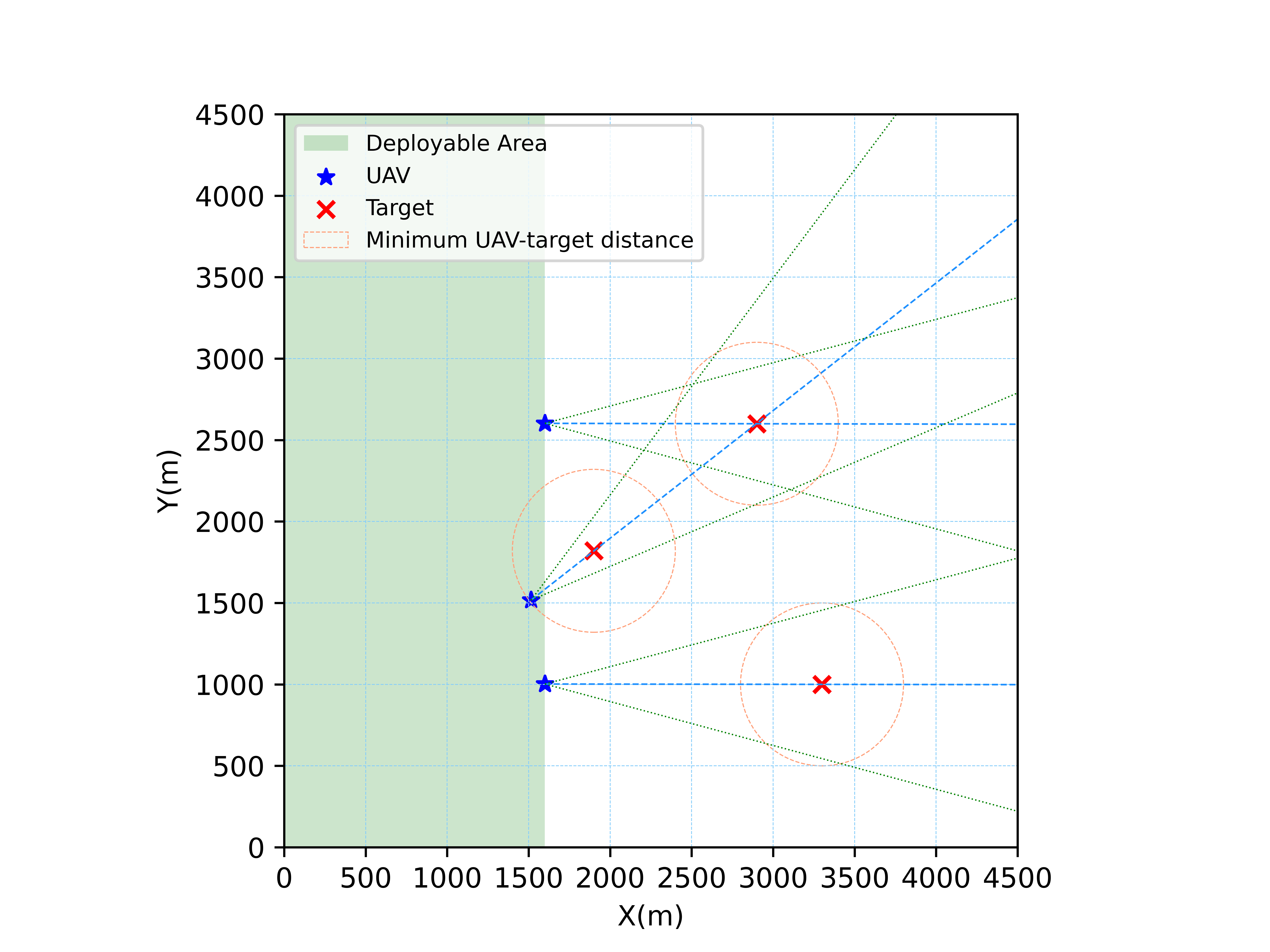}}
	\hspace{0.01in}
	\subfigure[$\text{ 4 UAVs cooperatively interfere with 3 targets}.$]{
		\label{fig:subfig:d}
		\includegraphics[width=8.8cm,height=2.5in]{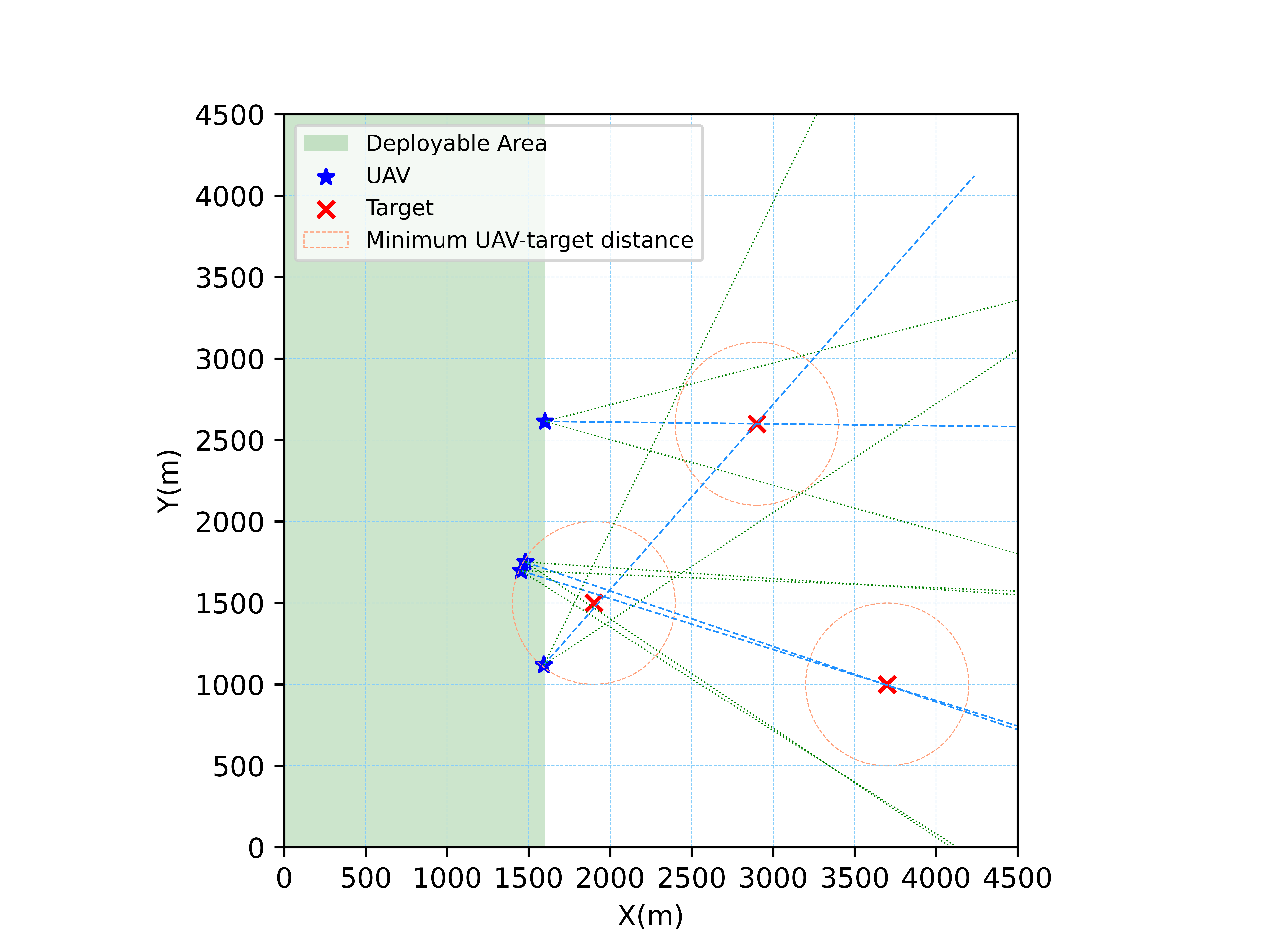}}
	\hspace{0.01in}
	\caption{Deployment and antenna orientation results in different scenarios.}
	\label{fig:UAV_deployment}
\end{figure*}

This section presents numerical simulations to verify the effectiveness of our proposed algorithm.
The scenario consists of $M$ UAVs to perform jamming tasks against $K$ targets. Each UAV is equipped with a directional antenna with a fixed beamwidth of $2\theta=30$ degrees, deployed at the same altitude $H = 600 \mathrm{~m}$. 
For simplification, the horizontal deployable area $D_q$ of the interference UAVs is modeled as a 2D space with $x \leq 1600 \mathrm{~m}$. 
Additionally, the minimum distance between the targets and the UAVs is set at $S_l = 500 \mathrm{~m}$, while the anti-collision distance of the UAVs is set at $R_l = 50 \mathrm{~m}$. 
 Other parameter settings are as follows: $N_t=5, \beta_0=-30 \mathrm{~dB}, \alpha^2=-110 \mathrm{~dBm}, P = 4\times 10^{-3} \mathrm{~W}, P_s = 2\times 10^{-2} \mathrm{~W}$. Initial values concerned with ADMM are set as follows: $\varpi_\chi = \varpi_\mu =200$, $\rho_1 = \rho_2 = 0.01$, $\eta =10^{{-}3}$.

Fig.\ref{fig:UAV_deployment} displays strategic interference UAV deployment aligned with suitable antenna orientations to optimize target interference. To demonstrate the efficacy of the proposed algorithm, we formulate scenarios where the targets are put at different positions and the control center is settled accordingly. In (c) and (d) of the Fig.\ref{fig:UAV_deployment}, when a UAV interferes with single target, it occupies an edge position of the deployable region, mirroring the target's y-coordinate, while its directional antenna aligns parallel to the x-axis, achieving optimal interference. In (b) and (d), UAVs go as near as possible to the targets while having their antenna main lobes covering the targets meticulously, to underpin the maximization of interference efficiency. Notably, when multiple UAVs converge to collectively target a shared objective, they exercise careful coordination to maintain a minimal inter-UAV distance while synchronizing their antenna orientations, as shown in (d). In scenarios like (a), in order to have the antenna main lobe cover more targets, UAV sacrifice distance to gain overall better interference. This strategic orchestration effectively magnifies collaborative interference, ultimately augmenting the impact on the targeted objective.

In order to evaluate the performance of the proposed algorithm (referred to as the "Proposed Scheme") for interference resource scheduling in scenarios involving multiple UAVs, a comparison is made against two baseline schemes:

• Baseline Scheme 1: This algorithm deploys each interference UAV in the deployable area, nearest to a target, while ensuring that the antenna mainlobe of each UAV covers as more targets as possible. Subsequently, a block coordinate descent algorithm is employed to optimize the antenna directions of the interference UAVs \cite{uav_wcnc}.

• Baseline Scheme 2: In this algorithm, deployment and antenna orientation is alternately updated to minimize the target function. Block coordinate descent algorithm is employed to optimize both the antenna directions and UAV positions.

\begin{figure}[!h]
	\centering
	\includegraphics[width=8.0cm,height=2.5in]{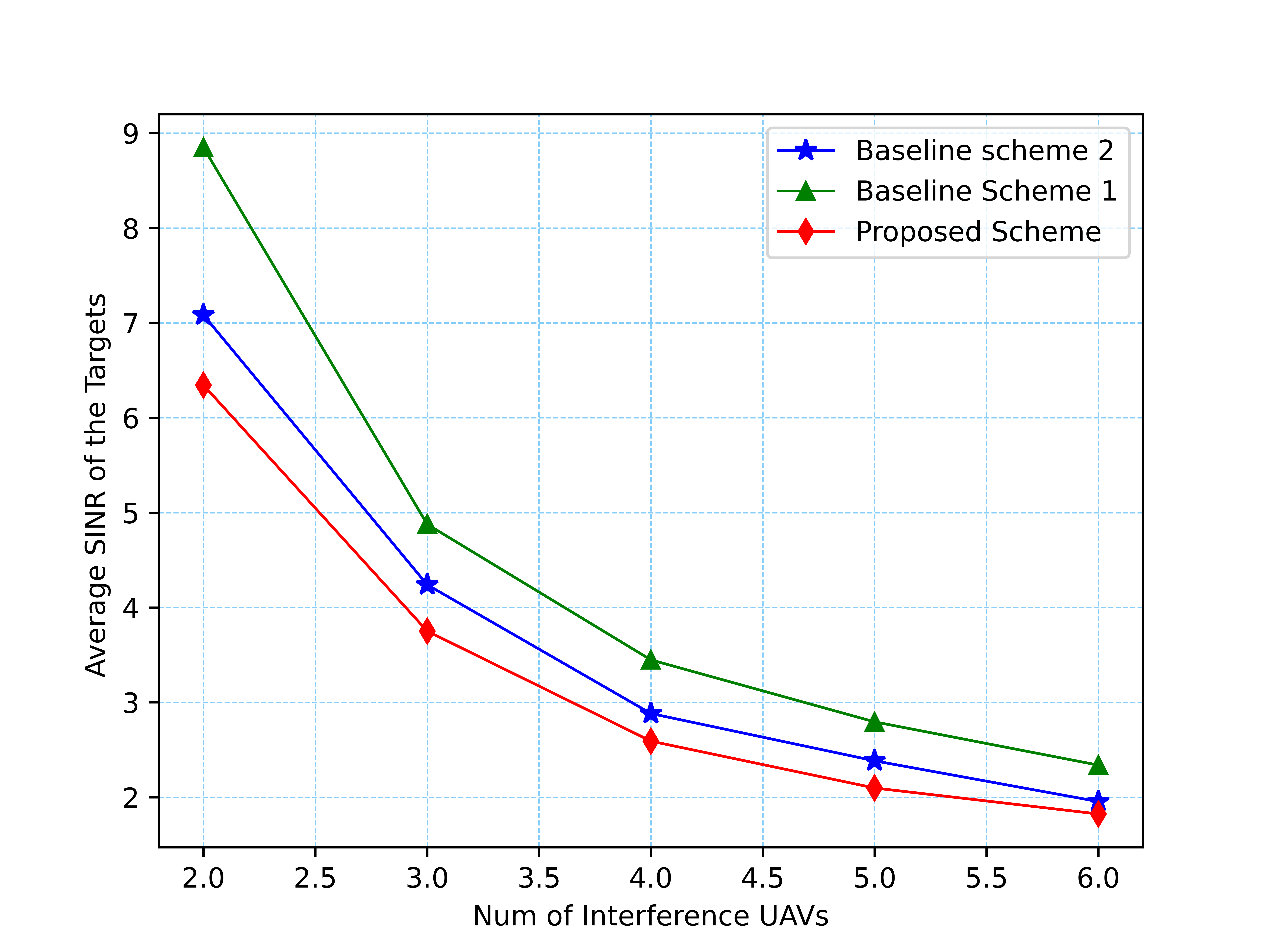}
	\caption{Average SINR of the targets versus the number of interference UAVs $M$.} 
	\label{fig:UAV_Performance}
\end{figure}

Experimental outcomes in Fig.\ref{fig:UAV_Performance} illustrate that as the number of interference UAVs increases, the average SINR decreases significantly for all algorithms, confirming the effectiveness of collaborative interference by multiple UAVs, while indicating that the proposed algorithm outperforms the other two baseline schemes.
This superiority arises from the algorithm presented herein, which undertakes a comprehensive evaluation of the collective influence exerted by the placements of interference UAVs and the orientations of antennas. This approach enhances the interference efficacy of the multi-UAV network.

\section{CONCLUSION}
In this paper, we jointly optimize the UAV deployment and directional antenna orientation in a UAV swarm-based interference network so as to minimize the communication performance of the targets. 
To tackle the formulated non-convex problem, we propose an efficient iterative algorithm for the solution with the aid of the ADMM, which decomposes the formulated non-convex problem into multiple subproblems and solves them alternately. And to accelerate convergence, we combined two machine learning optimizer NAG and RMSProp with the gradient projection method.
Extensive simulation results demonstrate the superior interference performance achieved by our proposed scheme, comparing with other benchmark methods.

\end{document}